\documentclass[a4paper, 11pt,reqno]{amsart}

\usepackage{amsmath,amsthm,amsbsy,amssymb}
\usepackage{arydshln} % dashline in table
\usepackage{array}
\usepackage{booktabs} %\toprule, \midrule \addlinespace
\usepackage{lipsum}
\usepackage{comment}
\usepackage{listings}

\makeatletter
\newcommand{\leqnos}{\tagsleft@true\let\veqno\@@leqno}
\newcommand{\reqnos}{\tagsleft@false\let\veqno\@@eqno}
\reqnos
\makeatother

\usepackage{xcolor} 
\definecolor{orange}{rgb}{1,0.5,0}
\definecolor{Red}{rgb}{.795,0.015,0.017}
\definecolor{Ggreen}{rgb}{0.,0.675,0.0128}
\definecolor{Bblue}{rgb}{0.16,.32,0.91}
 \usepackage[backref=page]{hyperref}
\hypersetup{
    colorlinks = true,
linkcolor={red},
urlcolor={blue},
citecolor={Ggreen},    
urlcolor = {blue},
citebordercolor = {0.33 .58 0.33},
 linkbordercolor = {0.99 .28 0.23},
 pdftitle={More than 8/9 zeros are simple or on the critical line},
 pdfauthor={CZ}
}

\usepackage{mathrsfs}% scriptfont

\newcommand{\cC}{\mathcal C}
\newcommand{\cD}{\mathcal D}

\newcommand{\cR}{\mathcal R}

\usepackage{bm}
\newcommand*{\B}[1]{\ifmmode\bm{#1}\else\textbf{#1}\fi}

\newcommand{\bg}{\B{g}}
\newcommand{\br}{\B{r}}
\newcommand{\bc}{\B{c}}

\newcommand{\RR}{\mathbb{R}}
\newcommand{\ZZ}{\mathbb{Z}}

\newcommand{\HS}{\mathrm{HS}}
\newcommand{\cMT}{C_{\mathrm{MT}}}
\newcommand{\un}{N_{s\cup 0}}

\usepackage{mathtools}

\usepackage{stmaryrd}
\newcommand{\sdfrac}[2]{\mbox{\small$\displaystyle\frac{#1}{#2}$}}

\usepackage{xspace}

\newcommand*\dd{\mathop{}\!\mathrm{d}}

\usepackage{enumitem} % adjust spacings
\makeatletter
\def\namedlabel#1#2{\begingroup
    #2%
    \def\@currentlabel{#2}%
    \phantomsection\label{#1}\endgroup
}
\makeatother

\usepackage{float}
\newcommand{\tr}{\operatorname{tr}}

\numberwithin{equation}{section}

\theoremstyle{plain}
\newtheorem{theorem}{Theorem}

\newtheorem{lemma}{Lemma}[section]

\theoremstyle{remark}
\newtheorem{remark}{Remark}[section]
\newtheorem*{remark*}{Remark}
\newtheorem*{example*}{Example}
\newtheorem*{examples*}{Examples}

\theoremstyle{definition}

\usepackage{subfigure}  % package for subfigure formatting
\usepackage{float}
\usepackage[pdftex]{graphicx}
\graphicspath{{IMAGES-DENOMINATORS}{IMAGES-ORBITS/}{IMAGES/}{./}}
\usepackage[]{caption}
\renewcommand*{\backref}[1]{}
\renewcommand*{\backrefalt}[4]{%
  \ifcase #1 %
No citations.% use \relax if you do not want the "No citations" message
  \or
(page #2).%
  \else
(pages #2).%
  \fi%
}
\usepackage{cite}

\usepackage{tikz,}
\usetikzlibrary{trees,positioning}
\usetikzlibrary{math} %needed tikz library
\usetikzlibrary {shapes.multipart}

\usepackage{nicematrix}
\makeatletter
\renewcommand*\env@matrix[1][\arraystretch]{%
  \edef\arraystretch{#1}%
  \hskip -\arraycolsep
  \let\@ifnextchar\new@ifnextchar
  \array{*\c@MaxMatrixCols c}}
\makeatother

\usepackage{stmaryrd}
\usepackage[]{mdframed}

\begin{document}
\title[More than 8/9 zeros are simple or on the critical line]
{More than 8/9 zeros of the Riemann zeta function\newline are simple or on the critical line}

\author[Cristian Cobeli]{Cristian Cobeli}
\address[Cristian Cobeli]{``Simion Stoilow'' Institute of Mathematics 
of the Romanian Academy,~21 Calea Griviței Street, P. O. Box 1-764, Bucharest 014700, Romania}
\email{cristian.cobeli@imar.ro}

\author{Alexandru Zaharescu}
\address[Alexandru Zaharescu]{Department of Mathematics, 
University of Illinois at Urbana-Champaign, Urbana, IL 61801, USA,
and 
``Simion Stoilow'' Institute of Mathematics of the Romanian Academy,~21 
Calea Griviței 
Street, P. O. Box 1-764, Bucharest 014700, Romania}
\email{zaharesc@illinois.edu}

\makeatletter
\@namedef{subjclassname@2020}{\textup{2020} Mathematics Subject Classification}
\makeatother
\subjclass[2020]{11M06, 11M26}
\thanks{Key words and phrases: Riemann zeta function, critical line, simple zeros, pair correlation}
% \keywords{Riemann zeta function, simple zeros, critical line, pair correlation, Gram matrix, interval arithmetic}

\begin{abstract}
We prove unconditionally that the lower limiting proportion of 
nontrivial zeros of the Riemann zeta function, which are simple or lie 
on the critical line, is greater than~$8/9$, with zeros counted 
with multiplicity. 
\end{abstract}
%%% --------------------------------------------------
\maketitle
% {
%   \hypersetup{linkcolor=blue}
% \hypersetup{pdftitle={More than 8/9 zeros of the Riemann zeta function 
% are simple or on the critical line},pdfauthor={Cristian Cobeli and Alexandru Zaharescu}}
%   \small
%   \tableofcontents
% }
%%% --------------------------------------------------

%%%%%%%%%%%%%%%%%%%%%%%%%%%%%%%%%%%%%%%%%%%%%%%%%%%%%%%%%%%%%%%
\section{Introduction}
Riemann's 1859 statement~\cite{Rie1859} 
on the non-trivial zeros
of the zeta function, according to which it is 
 ``\emph{sehr wahrscheinlich}''
% very probable 
that all lie on the critical line, 
became known as the \emph{Riemann hypothesis}.
Let $N(T)$ denote the number of  nontrivial zeros $\rho = \beta+i\gamma$ of 
$\zeta(s)$ with $0<\gamma\le T$, counted with multiplicity. 
Write $m_\rho$ for the multiplicity at a distinct zero location $\rho$. 
Let $N_0(T)$ count the zeros on the critical line with multiplicity, 
$N_s(T)$ count simple zeros, and $N_0^s(T)$ count simple zeros on the critical line. 
Thus
\begin{equation*}%\label{equation_uniondef}
   \un(T) = N_s(T)+N_0(T)-N_0^s(T)
\end{equation*}
counts zeros that are simple, or on the critical line, or both. 
A multiple critical zero contributes its full multiplicity to this union. Its complement is
\begin{equation}\label{equation_badmass}
   B(T) = \sum_{\substack{\rho\ \mathrm{distinct},\ 0<\gamma\le T\\
                      \beta\ne1/2,\ m_\rho\ge2}}m_\rho
                      \ \ \text{ and }\ \ 
 \un(T) = N(T)-B(T).
\end{equation}

Hardy~\cite{Har2014} proved in 1914 that 
infinitely many zeros of the zeta function lie on the critical line. 
Then Hardy and Littlewood~\cite{HL2021} strengthened this to $N_0(T)\gg T$. 
Selberg~\cite{Sel1943} introduced mollification into the problem 
and proved that a positive proportion of the zeros are on the critical line. 
Levinson~\cite{Lev1974} obtained a proportion greater than $1/3$, 
and Conrey~\cite{Con1989} raised it above $2/5$. 
Further developments of twisted moments and mollifiers include the works of 
Bui, Conrey and Young~\cite{BCY2011}, 
Feng~\cite{Fen2012},
Robles, Roy and the second author~\cite{RRZ2016},
Pratt and Robles~\cite{PR2018}.
Wu~\cite{Wu2019} and independently
Pratt, Robles, Zeindler, and the second author~\cite{PRZZ2020} subsequently proved 
that more than $5/12$ of the zeros are on the critical line. 

Montgomery~\cite{Mon1973}, 
assuming the Riemann hypothesis,
as an application of 
his pair correlation theorem, established a proportion of at least $2/3$ of simple zeros.
Then,
Montgomery and Taylor chose a test function 
(see Montgomery's remarks toward the end of~\cite{Mon1975}) that improved the proportion to
\begin{equation}\label{equation_MT}
   2-\cMT = 0.6725007\ldots,\ \text{ where }
   \ \cMT = \frac 12+\frac 1{\sqrt2}\cot\frac 1{\sqrt2}.
\end{equation}
Aryan~\cite{Ary2022} obtained an unconditional Fej\'er-kernel formula, 
and Baluyot, Goldston, Suriajaya and Turnage-Butterbaugh~\cite{BGSTB2024} 
established an unconditional form of Montgomery's result. 
Their subsequent work~\cite{BGSTB2026} connected the
method with horizontal multiplicity and obtained 
conclusions for zeros on the critical line 
under hypotheses on a narrow vertical box centered on the critical line. 
In particular, their Remark~2 records an $8/9$ union proportion under such a hypothesis. 
Goldston and Suriajaya~\cite{GS2026a,GS2026b} further explained this mechanism. 
Goldston, Lee, Schettler, and Suriajaya~\cite{GLSS2026} 
showed that the full pair correlation conjecture, 
without separately assuming the Riemann hypothesis, 
implies that asymptotically almost all zeros are simple and lie on the critical line.

A proof discovered and written by Claude~\cite{Ant2026}, 
and verified by Alp\"oge and Furman~\cite{AF2026}, 
showed that the proportion of zeros that are simple and lie on the critical line
is larger than or equal to the quantity in~\eqref{equation_MT} unconditionally. 
Lamzouri~\cite{Lam2026} gave a simpler proof for the same result using 
a finite Hilbert space inequality and the unconditional pair correlation formula. 
He also proved the union lower bound
\begin{equation}\label{equation_Lamzouri}
   \liminf_{T\to\infty}\frac{\un(T)}{N(T)} 
   \ge  1-\frac{\cMT-1}{3/2+\sqrt2}  =  0.88762\dots.
\end{equation}
We establish the following unconditional strengthening.

%%%%%%%%%%%%%%%%%%%%%%%%%%%%%%%%%%
\begin{theorem}\label{Theorem_main}
The nontrivial zeros of the Riemann zeta function satisfy
\begin{equation*}%\label{equation_main}
   \liminf_{T\to\infty}\frac{\un(T)}{N(T)}>\frac89.
\end{equation*}
Consequently, $\un(T)>(8/9)N(T)$ for every sufficiently large $T$.
\end{theorem}

The proof of Theorem~\ref{Theorem_main} produces the constant in~\eqref{equation_U}, which is 
approximately $0.88893$. 
% 0.888931016871589 = U
%%%%%%%%%%%%%%%%%%%%%%%%%%%%%%%%%%%%%%%%%%%%%%%%%%%%%%%%%%%%%%%
\section{Lamzouri's proof in brief}\label{Section_Lamzouri}
For a real, even, compactly supported function $\eta$, normalized by $\int\eta^2 = 1$, put
\begin{align}\label{equation_notation}
   F = \eta^2,\qquad K(z) = \widehat F(z) = \int_\RR F(u)e^{-2\pi iuz}\dd u.
\end{align}
Let $Z$ be a finite multiset invariant under conjugation, with $m_z = m_{\bar z}$. Unless otherwise indicated, sums indexed by $z\in Z$ below are over distinct locations, 
and the multiplicities are written explicitly. Define
\begin{equation}\label{equation_energy}
   N = \sum_zm_z\ \ \text{ and }\ \ E_K(Z) = \sum_{z,w}m_zm_wK(z-w)^2.
\end{equation}
With $f_z(u) = \eta(u)e^{-2\pi iuz}$, conjugation symmetry gives
\begin{equation}\label{equation_tensor}
   E_K(Z) = \Big\|\sum_zm_z f_z\otimes f_z\Big\|_{L^2(\RR^2)}^2.
\end{equation}
Lamzouri groups the vectors belonging to real multiple locations 
and nonreal conjugate pairs, extends the resulting nested subspaces 
by the simple real vectors, and applies Bessel's inequality in an adapted orthonormal basis. 
The square inequalities used on the successive subspaces give
\begin{equation}\label{equation_Lamfinite}
   n\ge2N-E_K(Z)\ \ \text{ and }\ \  N_d(Z)\ge\frac{3N-E_K(Z)}2,
\end{equation}
where $n$ is the number of simple real elements and $N_d(Z)$ the number of distinct elements 
(see~\cite[Proposition~2.1]{Lam2026}). 
The same Hilbert space framework was then used to prove~\eqref{equation_Lamzouri}.

For zeta zeros, take
\begin{equation}\label{equation_ZT}
   Z_T = \left\{\frac{i(\rho-1/2)\log T}{2\pi}:0<\gamma\le T\right\},
\end{equation}
with multiplicity. Reflection in the critical line becomes conjugation, 
and critical zeros become real elements. For each fixed smooth $\eta$ 
supported in $(-1/2,1/2)$, the unconditional pair correlation formula gives
\begin{equation}\label{equation_momentlimit}
 \frac{E_K(Z_T)}{N(T)} = \cC(f)+o_\eta(1)\ \ \text{ and }\ \ 
 \cC(f) = \int f(u)^2\dd u +
                 \iint|u-v|f(u)f(v)\dd u\dd v.
\end{equation}
There is a rational weight in the original formula. 
Lamzouri removes it by applying that formula separately to the two fixed 
functions $f*f$ and $(f*f)''$, as recalled in Section~\ref{Section_completion}. 
Smooth approximation of the Montgomery--Taylor density makes $\cC(f)$ tend to $\cMT$. 
Substitution in~\eqref{equation_Lamfinite} 
proves the simple-critical and distinct-zero bounds.

%%%%%%%%%%%%%%%%%%%%%%%%%%%%%%%%%%%%%%%%%%%%%%%%%%%%%%%%%%%%%%%
\section{The refinement leading past eight ninths}\label{Section_ideas}
We use only the pair correlation formula and the standard 
zero-counting asymptotic as external analytic inputs. 
There are three ingredients that allow us to pass over the $8/9$ barrier. 
\begin{itemize}[wide, labelwidth=0pt]
   \setlength{\itemsep}{0pt}
   \itemindent = 4pt
    \setlength{\parskip}{7pt}
    \setlength{\parsep}{0pt} 
    \item[\texttt{(A)}]\hypertarget{PartA}{}
First, we retain 
a spectral remainder of the simple-real Gram matrix $H = (K(x_i-x_j))$:
\begin{align*}
   \cR_\tau(H) = \tr\psi_\tau(H)\ \ \text{ and }\ \
   \psi_\tau(t) = (t-1)^2-
   (\max(0,(t-\tau)))^2
\ \text{ for }\ \tau\ge 2.
\end{align*}
% It is nonnegative, but it contains information lost by replacing it by zero. 
The inequality proved in Section~\ref{Section_spectral} below is
\begin{equation}\label{equation_roadmap}
   E_K(Z)\ge A_\tau N+(1-A_\tau)n + \frac{\tau^2}{4}B 
      + \cR_\tau(H)\ \ \text{ and }\ \ A_\tau
   = 2\tau-\frac{\tau^2}{2},
\end{equation}
where $B$ is the mass of multiple nonreal elements. 
Dropping the nonnegative 
remainder~$\cR_\tau(H)$
and taking $\tau = 2+\sqrt2$ 
gives the union coefficient in~\eqref{equation_Lamzouri}.
    \item[\texttt{(B)}]\hypertarget{PartB}{}
Second, an explicit inequality for seven ordered real points 
bounds a weighted sum of their squared correlations 
plus a small multiple of their total span (see Lemma~\ref{Lemma_seven} below). Then in Section~\ref{Subsection_Blocks}, summing this inequality 
inside consecutive blocks of $120$ points, and averaging 
the $120$ disjoint-block partitions, gives
\begin{equation}\label{equation_roadmapR}
   \cR_\tau(H)\ge\alpha n-dL-119\alpha,\ \ \text{ where }\ \
 \alpha = \frac{209}{5000},\ \ d = \frac{57}{2000},
\end{equation}
when the points lie in an interval of length $L$ and $\tau$ 
has the value below. Only scalar convexity and the spectral theorem 
are needed for this passage from local correlations to the full matrix.
    \item[\texttt{(C)}]\hypertarget{PartC}{}
Third, we choose
\begin{equation}\label{equation_tau}
 \tau=2+\sqrt{2-2\alpha} = 2+\frac{\sqrt{4791}}{50}.
\end{equation}
Then $A_\tau = 1+\alpha$, so the coefficient of $n$ cancels in~\eqref{equation_roadmap}. 
The resulting inequality bounds~$B$ directly. 
A three-term cosine perturbation of the Montgomery--Taylor density 
makes the geometric gain exceed its small increase in scalar moment. 
We will improve the combined moment-and-remainder bound, 
not the scalar moment alone.
\end{itemize}

%%%%%%%%%%%%%%%%%%%%%%%%%%%%%%%%%%%%%%%%%%%%%%%%%%%%%%%%%%%%%%%
\section{A finite spectral inequality}\label{Section_spectral}
We work with~\eqref{equation_energy}. Let $r$ be the number of multiple real locations, 
$k$ the number of nonreal conjugate pairs, and $b = r+k$. 
Since a multiple nonreal pair uses at least four units of mass, 
whereas every other block outside the $n$ simple real locations uses at least two,
\begin{equation}\label{equation_blockcount}
   b\le\frac{N-n}{2}-\frac B4.
\end{equation}
Indeed, if $q_B$ denotes the number of multiple nonreal pairs, then
$b\le(N-n-B)/2+q_B$ and $q_B\le B/4$.

Here is the exact operator underlying~\eqref{equation_tensor}. Put
\begin{align*}
   g_z = \frac{f_z+f_{\bar z}}2,\qquad
   h_z = \frac{f_z-f_{\bar z}}{2i}.
\end{align*}
The vectors $f_x$ for real $x$, and $g_z,h_z$, 
lie in the real Hilbert subspace of complex $L^2(\RR)$ 
consisting of functions $v$ with $v(u) = \overline{v(-u)}$. 
Their usual complex inner products are real. Direct integration gives
\begin{align*}
   \|f_x\|^2 = 1,\qquad 
   \|g_z\|^2-\|h_z\|^2 = 1,\qquad 
   \langle f_x,f_y\rangle  =  K(x-y).
\end{align*}
On their finite real span define the self-adjoint operator
\begin{equation*}%\label{equation_operator}
   A = \sum_{x\in Z\cap\RR}m_x f_x\otimes f_x +
       2\sum_{\{z,\bar z\}}m_z(g_z\otimes g_z-h_z\otimes h_z).
\end{equation*}
Here $v\otimes v$ denotes the real rank-one operator 
$w\mapsto\langle w,v\rangle v$.
The trace of $A$ is $N$ and its tensor is the one in~\eqref{equation_tensor}.
Then the expansion gives
\begin{align*}
   \|A\|_{\HS}^2 = \sum_{z,w}m_zm_wK(z-\bar w)^2 = E_K(Z),
\end{align*}
where the last equality reindexes by conjugation.

Let $P$ be the contribution of the simple real locations and $Q = A-P$. 
Then $P\succeq0$, $\tr P = n$, $\tr Q = N-n$, and $n_+(Q)\le b$, 
where $n_+$ counts positive eigenvalues. 
The last assertion follows because each remaining real block 
has positive index one and each nonreal pair has positive index at most one: 
their sum is the pullback of their direct sum. 
The nonzero eigenvalues of $P$ agree with those of the Gram matrix $H$.

%%%%%%%%%%%%%%%%%%%%%%%%%%%%%%%%%%%%
\begin{lemma}\label{Lemma_square}
For self-adjoint $P\succeq 0$ and $Q$ with $n_+(Q)\le b$, and for every $\tau\ge2$,
\begin{equation*}%\label{equation_square}
   \|P+Q\|_{\HS}^2 \ge 2\tau\tr Q-\tau^2b+\tr\phi_\tau(P),
\end{equation*}
where $\phi_\tau(t) = t^2 - (\max(0,(t-\tau)))^2$.
\end{lemma}
\begin{proof}
Let us write $Q = Q_+-Q_-$ with $Q_\pm\succeq0$ and $Q_+Q_- = 0$. Since $\tr(PQ_+) \ge 0$,
\begin{align*}
   \|P+Q\|_{\HS}^2 \ge \|P-Q_-\|_{\HS}^2+\|Q_+\|_{\HS}^2.
\end{align*}
Applying $q^2 \ge 2\tau q-\tau^2$ to the at most $b$ positive eigenvalues reduces the claim to
\begin{align*}
   \|P-Q_-\|_{\HS}^2+2\tau\tr Q_- \ge \tr\phi_\tau(P).
\end{align*}
Then we diagonalize $P$, discard the squared off-diagonal entries of $Q_-$, and use
\begin{align*}
   (p-x)^2+2\tau x \ge  p^2 - (\max(0,(p-\tau)))^2\
   \text{ for } p,x \ge 0.
\end{align*}
% The scalar minimum occurs at $x = \pos{p-\tau}$.
\end{proof}

Since $\phi_\tau(0) = 0$, $\tr H = n$, and $P,H$ have the same nonzero eigenvalues,
\begin{align*}
   \tr\phi_\tau(P) = n+\cR_\tau(H).
\end{align*}
Apply Lemma~\ref{Lemma_square} and then~\eqref{equation_blockcount}. 
This proves~\eqref{equation_roadmap}, as aimed in 
part~(\hyperlink{PartA}{\texttt{A}}). The function $\psi_\tau$ 
is nonnegative and convex: its derivative increases from $-2$ to $2(\tau-1)$ 
and is then constant. These elementary properties will suffice to bound the remainder.

%%%%%%%%%%%%%%%%%%%%%%%%%%%%%%%%%%%%%%%%%%%%%%%%%%%%%%%%%%%%%%%
\section{
An explicit kernel and its geometric remainder
}\label{Section_kernel}
\subsection{Kernel selection and numerical implications
% The kernel and the finite numerical assertion
}
\begin{figure}
\begin{minipage}[t]{0.49\linewidth}
    \centering
    \includegraphics[width=0.99\textwidth, height=0.47\textwidth]{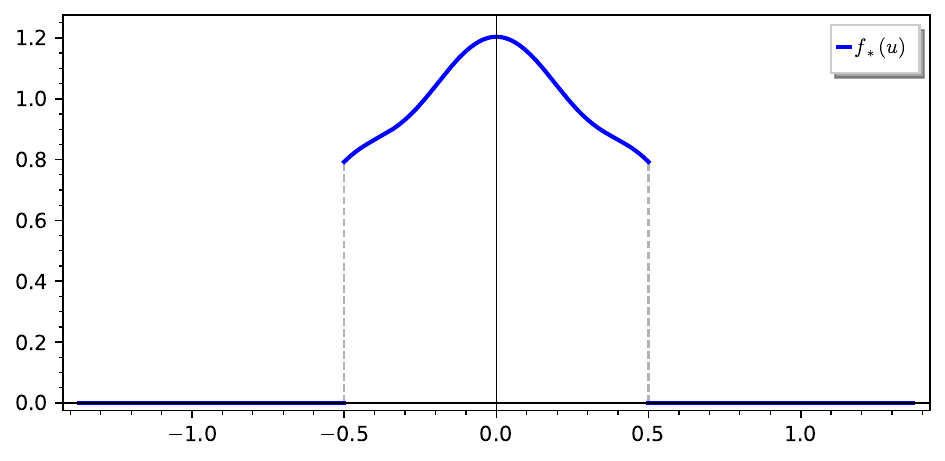}
    \captionsetup{width=.99\linewidth}
    \caption{The density kernel $f_*(u)$.
    }
    \label{Figure_kernelfstar}
\end{minipage}
\hfill
\begin{minipage}[t]{0.49\linewidth}
%%%%%%%%%%%%%%%%%%
    \centering
    \includegraphics[width=0.99\textwidth, height=0.47\textwidth]{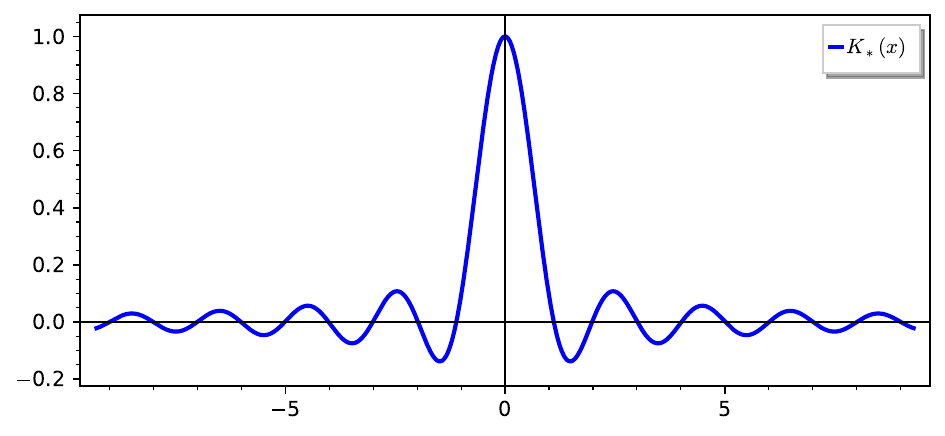}
    \captionsetup{width=.99\linewidth}
    \caption{\mbox{The Fourier transform $K_*(x) = \widehat f_*(u)$}.
    }
    \label{Figure_FTKstar}
\end{minipage}
\end{figure}

We let $f_*(u)$ be defined on the interval $I = [-1/2,1/2]$ by
\begin{equation}\label{equation_density}
   f_*(u) = 
   \frac{\cos(\sqrt2u)}{\sqrt2\sin(1/\sqrt2)}
   +\frac2{25}\cos(2\pi u)+\frac 1{25}\cos(4\pi u)
   -\frac 1{200}\cos(6\pi u),
\end{equation}
and $f_*(u) = 0$ if $u\not\in I$.
The kernel 
$f_*(u) = h_0(u)+h_1(u)+h_2(u)+h_3(u)$ is structured as a primary base kernel 
$h_0(u) = \cos(\sqrt{2}u)/(\sqrt{2}\sin(1/\sqrt{2}))$,
which provides the core frequency confinement 
on $I$, combined with a trigonometric polynomial
with three terms 
\mbox{$h_k(u)=a_k \cos(2\pi k u)$}.
The coefficients 
$a_1 = 2/25$, 
$a_2 = 1/25$, 
and $a_3 = -1/200$ 
are weights chosen to tune the shape of the Fourier transform 
$K_*(x) = \widehat f_*(u)$.
As shown in the graphical representation of $f_*(u)$ in Figure~\ref{Figure_kernelfstar},
the minimal value of $f_*(u)$ on $I$ is 
\mbox{$f_*(\pm 1/2) = (1/\sqrt{2})\cot(1/\sqrt{2})-7/200
% \approx 0.79249
$}, which is 
larger than $0.79$.

The kernel satisfies the requirement 
$\int f_* = \int_I f_0 = 1$.
This means that the normalization condition  is automatically 
satisfied independently of the choice of the coefficients $a_1,a_2,a_3$. 
These coefficients are chosen to ensure local-correlation estimates
and the admissible approximations introduced in Section~\ref{Section_completion}.

%%%%%%%%%%%%%%%%%%%%%%%%%%%%%%%%%%%%%%%%%%
Since $f_*(u)$ is an even function, 
its Fourier transform $K_*(x) = \widehat{f}(x)$
is real for $x$ real. Thus, while integrating,
the sine terms vanish, leaving us with the following expression:
\begin{align*}
   K_*(x) &=   \int_{-1/2}^{1/2} f(u)\cos(2\pi ux) \dd u\\
   &=
   \sdfrac{1}{\sqrt2 \sin(1/\sqrt2)} \left(
   \sdfrac{\sin((\sqrt2 - 2\pi x)/2)}{\sqrt2 - 2\pi x} + \sdfrac{\sin((\sqrt2 + 2\pi x)/2)}{\sqrt2 + 2\pi x} \right)+
   \sdfrac{x\sin(\pi x)}{\pi}\sum_{k=1}^3\sdfrac{(-1)^ka_k}{x^2-k^2}.
\end{align*}
For $0\le x \le 1$, differentiation under the integral gives
\begin{equation*}
   K_*'(x) = -4\pi\int_0^{1/2}u f_*(u)\sin(2\pi ux)\dd u \le 0,
\end{equation*}
so that $K_*(x)$ is decreasing
 on $[0, 1]$ (see the graphic representation in Figure~\ref{Figure_FTKstar}). 
The exact value of $K_*(1/2)$ is 
\begin{align}\label{equation_Kstar12}
   K_*(1/2) = \frac{2\pi(\cMT-1/2)}{\pi^2-2}+\frac{167}{3500\pi}
   \approx 0.67587.
   % \approx 0.675873162235178. % > 2/3.
\end{align}

%%%%%%%%%%%%%%%%%%%%%%%%%%%%%%%%%%%%%%%%%%%%%%
%%%%%%%%%%%%%%%%%%%%%%%%%%%%%%%%%%%%%%%%%%%%%%%%

%%%%%%%%%%%%%%%%%%%%%%%%%%%%%%%%%%%%%%%%%%%%%%
\begin{lemma}\label{LemmaMomentC}
With the kernel density $f_*$ defined by~\eqref{equation_density},
the value of the moment from~\eqref{equation_momentlimit}~is
\begin{equation}\label{equation_C}
  \begin{split}
   C: = \cC(f_*) 
    &= \int f_*(u)^2\dd u +
                 \iint|u-v|f_*(u)f_*(v)\dd u\dd v \\ 
   &= \cMT+\frac{321}{80000}-\frac{2449}{1440000\pi^2}.
  \end{split}
\end{equation}
\end{lemma}
\begin{proof}
We write $C$ as
\begin{align*}
   C= \int_{-1/2}^{1/2} f_*(u)^2 \dd u 
   + \langle f_*, \Psi f_* \rangle,
\end{align*}
where $\Psi $ is the integral operator that acts on functions supported on $I$ by
\begin{align*}
   (\Psi f)(u) = \int_{-1/2}^{1/2} |u-v| f(v) \dd v. 
\end{align*}
For $\omega\neq 0$, the action of~$\Psi$ on any harmonic $\cos(\omega u)$ is given by 
\begin{align}\label{equation_Psi}
   \Psi(\cos(2\pi k\cdot))(u) 
   \frac{\sin(\omega/2)}{\omega}
   +\frac{2}{\omega^2} \left( \cos\left(\frac{\omega}{2}\right) - \cos(\omega u) \right). 
\end{align}
This can be verified by differentiating it twice and matching the value and derivative at~\mbox{$u = 1/2$}.
For $u\in [-1/2,1/2]$, we write 
$f_*(u) = h_0(u) + h(u)$, where $h(u)=h_1(u) + h_2(u) + h_3(u)$, 
\begin{align*}
   h_0(u) &= \frac{\cos(\sqrt{2}u)}{\sqrt{2}\sin(1/\sqrt{2})},\ \ \ \ 
   h_k(u)=a_k \cos(2\pi k u) \ \ \text{ for $k=1,2,3$,}
   % h_1(u) = \frac{2}{25}\cos(2\pi u),\ 
   % h_2(u) = \frac{1}{25}\cos(4\pi u),\ 
   % h_3(u) = -\frac{1}{200}\cos(6\pi u).
\end{align*}
and $a_1 = 2/25$, 
$a_2 = 1/25$, 
 $a_3 = -1/200$.
Then $\int_I h_0=1$ and $\int_I h=0$. Since $(\Psi h_0)''=2h_0$ and 
$h_0''=-2h_0$, the function $(I+\Psi)h_0$ is affine and even, hence constant. 
Moreover
\begin{align*}
   (\Psi h_0)(1/2)=\frac12\int_Ih_0-\int_Iv h_0(v)\dd v = \frac12,
\end{align*}
so that
\begin{equation*}
   (I+\Psi)h_0=h_0(1/2)+\frac12=\cMT.
\end{equation*}
The operator is self-adjoint because its kernel $|u-v|$ is symmetric. Thus
\begin{align*}
    \langle h_0,(I+\Psi)h_0\rangle=\cMT,\ \ \text{ and }
   \ \ \langle h,(I+\Psi)h_0\rangle=\cMT\int_Ih=0.
\end{align*}
For $k\ge 1$, formula~\eqref{equation_Psi} specializes to
\begin{align*}
 \Psi(\cos(2\pi k\cdot))(u)
   =\frac{(-1)^k-\cos(2\pi ku)}{2\pi^2k^2}.
\end{align*}
Different integer-frequency cosines are orthogonal, have mean zero, 
and have squared integral $1/2$. Therefore all their cross terms also vanish, giving
\begin{align*}
   C(f_*) &=\langle f_*,(I+\Psi)f_*\rangle\notag\\
          &=\cMT+\frac12\sum_{k=1}^3 a_k^2
            \left(1-\frac{1}{2\pi^2k^2}\right)\notag\\
        &= \frac 12+\frac 1{\sqrt2}\cot\frac 1{\sqrt2} +
       \frac{2\pi^2 - 1}{625\pi^2} +
       \frac{8\pi^2 - 1}{10000\pi^2} +
       \frac{18\pi^2 - 1}{1440000\pi^2},
\end{align*}
which reduces to~\eqref{equation_C}.
\end{proof}
%%%%%%%%%%%%%%%%%%%%%%%%%%%%%%
\begin{remark}\label{Remark_CMT-C}
A comparison of the approximate numerical values of the two constants gives
\begin{align*}
    \cMT = 1.3274992\dots\ \ \text{ and }\ \ C = 1.3313394\dots
    % \cMT\approx 1.32749929632059\ \ \text{ and }\ \ C\approx 1.33133947994646.
\end{align*}
   % $\cMT\approx 1.32749929632059$ and $C\approx 1.33133947994646$.
% CMT = 1.32749929632059
% C = 1.33133947994646

% 0.888888888888889 = 8/9
% 0.888931016871589 = U
\end{remark}

Next, to carry out part~(\hyperlink{PartB}{\texttt{B}}) of our plan, 
we need an inequality that estimates the lower bound of a weighted 
sum $S(\bg)$ defined in~\eqref{equation_seven} on the gaps 
between points in a six-dimensional space.
In the first part of the sum, $q(x):=(K_*(x))^2$ acts like a measuring tool 
for the size of each gap, while the second part ponders  $S(\bg)$ with a fraction of the largest gap.
The inequality we obtain is very tight, with a difference of less than 
three thousandths from a numerical estimate of the absolute minimum of the sum.
 % Minimum S  = 0.0469769074
 %  Best g     = [1.87680678 1.15185995 1.8793085  1.13039575 1.82542535 1.13620368]
 %  sum(g)     = 9.000000
 %  Target     = 0.0440000000
 %  Margin     = 2.9769074104e-03
 
In the following lemma, we will see that the analysis in a 6-dimensional 
continuum space can be reduced to a number of discrete evaluations at 
the corners of certain boxes. This can be done in two ways: one is to 
proceed with a direct analytical proof, and the other is to leave the 
laborious calculations to a computer for verification.
The first option, however, is technical and tedious, 
involving numerous approximations using Taylor polynomials.
The second requires some time and computing power because, 
for instance, to carry on the smaller values of $K_*(x)^2$,
a mesh division of a $6$-dimensional box of size $50$ would require
the evaluation of the sum in $15\,625$ million points.
%15625000000
We will take a middle approach, reducing the set to be analyzed 
as much as possible, showing the necessary evaluations for the 
higher-order derivatives of $q(x)$, and leaving the verification of the remaining cases to a 
lightweight verification that we completed in~\texttt{SageMath}.

Next, to define $S(\bg)$, we introduce the nonnegative weights 
$w_{ij}$, $0\le i<j\le 6$ with the precise values listed in
Table~\ref{Table_weights}.
The weights are chosen such that every diagonal sums to~$2$, that is,
\begin{equation}\label{equation_wdiagonalsum}
   \sum_{i = 0}^{6-j}w_{i,i+j} = 2\quad(1\le j\le6)
\end{equation}
and the total weight is $12$.

%%%%%%%%%%%%%%%%%%%%%%%%%%%%%%%%%%%%%%%%%%%%%%%%%%
\begin{center}
\begin{table}
\caption{The table of weights $w_{i,j}$, $0\le i\le 5$, $1\le j\le 6$.}
\label{Table_weights}
\begin{tabular}{@{}ccccccc@{}}
\toprule
$i$ & $w_{i,1}$ & $w_{i,2}$ & $w_{i,3}$ & $w_{i,4}$ & $w_{i,5}$ & $w_{i,6}$\\
\midrule
0 & 0.306 &  0.645 & 1.000 & 1.000 & 1.000 & 2.000\\
1 &       &  0.347 & 0.047 & 0.000 & 0.000 & 1.000\\
2 &       &        & 0.347 & 0.616 & 0.000 & 1.000\\
3 &       &        &       & 0.347 & 0.047 & 1.000\\
4 &       &        &       &       & 0.347 & 0.645\\
5 &       &        &       &       &       & 0.306\\
\bottomrule
\end{tabular}
\end{table}
\end{center}
%%%%%%%%%%%%%%%%%%%%%%%%%%%%%%%%%%%%%%%%%%%%%%%

%%%%%%%%%%%%%%%%%%%%%%%%%%%%%%%%%%%%%%%%%%%%%%%%%%%%%%%%%%%%%%%%%%%%%%%%%%%
\begin{lemma}\label{Lemma_seven}
For every $\bg=(g_1,\dots,g_6)\in\RR^6$, where
$g_1,\ldots,g_6 \ge 0$, set $y_0 = 0$, and 
$y_k = g_1+\cdots+g_k$ for $k=1,2,\dots,6$. Then
\begin{equation}\label{equation_seven}
   S(\bg):= \sum_{0\le i<j\le6}w_{ij}K_*(y_j-y_i)^2
           + \frac 1{200}\sum_{i=1}^6 g_i 
       \ge \frac{11}{250}.
\end{equation}
\end{lemma}
\begin{proof}
The integral representation
of $K_*(x)$ implies that $K_*'(x)\le0$ on $[0,1]$ 
(confirmed by the graphical representation in Figure~\ref{Figure_FTKstar}),
and from~\eqref{equation_Kstar12} we find 
that $K_*(1/2)>2/3$. If some $g_i\le1/2$, together with one adjacent pair term, 
contributes to the first sum of $S(\bg)$ at least
$(3/10)(2/3)^2 = 2/15>11/250$, because all the weights 
on the longest diagonal in Table~\ref{Table_weights} are larger than $3/10$. 
If $\sum_i g_i \ge 8.8$, 
the second term of $S(\bg)$
is at least $8.8/200 = 11/250$. 

It remains to check inequality~\eqref{equation_seven} 
for all tuples $\bg$ in $\Omega$, where 
\begin{equation*}%\label{equation_Omega}
   \Omega = \Big\{\bg\in[1/2,13/2]^6:\ \sum_i g_i\le 8.8\Big\}.
\end{equation*}

Further pruning can be done. Let $S_0=11/250=0.044$, and 
\mbox{$m=\min\limits_{0\le i\le5} w_{i,i+1}=0.306$.} 
If $S(g)<S_0$, then each adjacent-pair term must satisfy 
\begin{align*} 
   w_{i,i+1}K_*(g_{i+1})^2 < S_0. 
\end{align*} 
(Note that there are six adjacent pairs
$(0,1), \dots,(5,6)$, so an argument involving $w_{i-1,i}$ 
applies to every gap, including $g_6$.) 
Consequently, 
\begin{align*}%\label{equation_5555} 
   |K_*(g_{i+1})|< \sqrt{\sdfrac{S_0}{w_{i,i+1}}} \le \sqrt{\sdfrac{0.044}{0.306}} <0.3792.
\end{align*} 
  
However on $[1/2,1)$, $K_*$ is positive and decreasing, 
and a plain calculation shows 
\begin{align*} 
   K_*(0.743) > 0.38029 > 0.3792. 
\end{align*} 
Now, if any $g_k\le 0.743$, it follows that
\begin{align*}
   S(\bg)\ge 0.306 K_*(g_k)^2
      > 0.306\cdot 0.38^2 > 11/250.
\end{align*} 
% 0.38^2*0.306-11/250 = 0.000186400000000003
% 0.3702*0.3702*0.306-11/250 = -0.00206329976000000
% 0.3792*0.3792*0.306-11/250 = 5.47840000003630e-7
If $\sum g_k\ge44/5$, the linear term in $S(\bg)$ alone proves the conclusion.
Otherwise, every gap is greater than $0.743$, and each is less than
$44/5-5\cdot 0.743 = 5.085$. Thus, every possible counterexample 
belongs to a smaller region $\Omega'\subset\Omega$, where 
\begin{align*}%\label{equation_Omega1} 
   \Omega'= \Big\{\bg\in(0.743,5.085)^6: \sum_{i=1}^6g_i\le8.8 \Big\},
\end{align*} 
and this is the remaining region on which~\eqref{equation_seven}
still needs to be verified.

%%%%%%%%%%%%%%%%%%%%%%%%%%%%%%%%%%%%%%%%%%%%%%
Let $b=1/\sqrt{2}$, $s(z) = \sin(z)/z$ for $z\neq 0$, and $s(0)=1$. 
The following is a useful algebraic equivalent expression for $K_*(x)$:
\begin{align*}%\label{equation_Kstar2}
   K_*(x) = \frac{b}{2\sin b}
   \left(s(\pi x-b)+s(\pi x+b)\right)
    + \frac 12\sum_{k=1}^3
    a_k\left(s(\pi(x-k))+s(\pi( x+k)\right).
\end{align*}
This formula is valid for every $x\in\RR$ because the definition $s(0) = 1$ handles all the apparent
singularities.

We denote $q(x)=K_*(x)^2$.
Since, as we observed before, $f_*(u)>0.79$
and $\int f_* = 1$, it follows that~$K_*^2(x)$ has the representation
\begin{align*}
   q(x)=\int_{-1}^1 h(t)e^{2\pi ixt}\dd t,
\end{align*}
where  $h=f_* * f_*$,
is a nonnegative probability density supported on $[-1,1]$. 
Note that the sign in the exponential is irrelevant because 
both $h$ and $q$ are even functions, just as $K_*$ and $f_*$ were.
Differentiation under the integral gives
\begin{align*}
   q^{(r)}(x)= \int_{-1}^1(2\pi it)^r h(t)
   e^{2\pi ixt}\dd t\ \ 
   \text{ for $r \in\ZZ_{\ge 0}$}.
\end{align*}
Hence, for every real $x$,
\begin{align*}
   |q^{(r)}(x)|\le(2\pi)^r.
\end{align*}
In particular, we obtain the following useful bounds:
$|q(x)|\le 1$, $|q'(x)|\le 2\pi < 6.284$,  $|q''(x)|\le 4\pi^2 < 40$, 
$|q'''(x)| \le 8\pi^3 < 249$,  
$|q^{(4)}(x)| \le 16\pi^4 < 1600$ for all $x\in\RR$. 

It suffices to verify the inequality $S(\bg)\ge 11/250$
on a series of boxes that intersect and cover
the closed enclosure
$\cD=\left[0.743, 5.085\right]^6 \cap
   \big\{\bg:\sum_{i=1}^6g_i\le 44/5\big\}
   \supseteqq \Omega'.$

For each box $\prod[\ell_k,u_k]$, 
intersect with the sum constraint by replacing
\begin{align*}
   u_k \leftarrow 
   \min\Big(u_k,\sdfrac{44}{5}-\sum_{r\ne k}\ell_r\Big).
\end{align*}
All these endpoint contractions are repeated until no further endpoint changes.
We exclude boxes only when the intersection is empty. 
For a pair interval one may further replace its upper endpoint by
$\min(\sum_{k=i+1}^ju_k,44/5-\sum_{k\notin(i,j]}\ell_k)$.
A box is discarded if any interval becomes empty.
For each pair $i<j$, the gap $y_j-y_i$ belongs~to
\begin{align*}
   J_{ij}= \left[\sum_{k=i+1}^{j}\ell_k,\, \sum_{k=i+1}^{j}u_k \right].
\end{align*}
Intersect this interval with the restrictions implied by the sum constraint. 
In particular, any box for which
\begin{align*}
   \sum_{k=i+1}^{j}\ell_k+ \sum_{k\notin\{i+1,\ldots,j\}}\ell_k >\frac{44}{5}
\end{align*}
is empty.

In a continuous run, with the goal of ending up with no unchecked boxes 
on the list containing those where the inequality is not yet verified, 
we use two tests called test \texttt{L1(B)} and test \texttt{L2(B)}. 
In   test \texttt{L1(B)}, for a nonempty box $B$, one needs
to find a lower bound for
\begin{align*}
L_1(B) :=
\sum_{i<j}w_{ij}\inf_{x\in J_{ij}}q(x)
+\frac1{200}\sum_{i=1}^6\ell_i.
\end{align*}
This requires bounds for $q$ and $q'$ on each $J_{ij}$.
These can be obtained by proving monotonicity on 
iterated subintervals and evaluating the appropriate endpoints with outward rounding;.
If the resulting lower bound for $L_1(B)$ is 
at least $11/250$, then the inequality is proved on the box $B$.
If not, we move to the second test or split 
$B$ in two subboxes that need to be analyzed separately.

To describe test \texttt{L2(B)}, let a box $B$ 
have center $\bc=(c_1,\dots,c_6)$ and radii 
$\br=(r_1,\dots,r_6)$. 
Define
\begin{align}\label{equation_center-radii}
    s_{ij}=\sum_{k=i+1}^jc_k,
 \qquad r_{ij}=\sum_{k=i+1}^jr_k,
 \qquad
 D_k=\frac1{200}+\sum_{i<k\le j}w_{ij}q'(s_{ij}).
\end{align}
Taylor's theorem and $|q''|<40$ imply, for every point $\bg$ in this box,
\begin{equation}\label{equation_TaylorbBox}
   S(\bg)\ge S(\bc)-\sum_{k=1}^6|D_k|r_k
              -20\sum_{i<j}w_{ij}r_{ij}^2. 
\end{equation}
Then test \texttt{L2(B)}
evaluates $S(\bc)$ with a lower interval endpoint and 
$|D_k|$ with an upper endpoint, giving the right-side 
of~\eqref{equation_TaylorbBox}
as a lower bound. 

When neither test \texttt{L1(B)} nor \texttt{L2(B)} proves $S\ge 11/250$, 
bisect an interval and retain both children. 

A finite binary subdivision tree, 
with every nonempty box 
that has passed test \texttt{L1(B)} 
or~\texttt{L2(B)} provides a completion of the proof of 
the inequality~\eqref{equation_seven} on $\Omega'$, 
and thus concludes the proof of Lemma~\eqref{Lemma_seven}.
\end{proof}
A step by step description of the algorithm of 
the verification of~\eqref{equation_seven} on $\Omega'$ in given in
Appendix~\ref{Appendix}.

%%%%%%%%%%%%%%%%%%%%%%%%%%%%%%%%%%%%%%%%
\subsection{Averaging principal blocks}\label{Subsection_Blocks}
We now deduce~\eqref{equation_roadmapR}, retaining a little 
flexibility for smooth approximation. Suppose a real positive-definite 
kernel $K$, $K(0) = 1$, satisfies~\eqref{equation_seven} with $11/250$ 
replaced by $\epsilon>0$, and put $p = 1/200$. Let $H = (K(x_i-x_j))$ 
for increasing real points in an interval of length $L$.

For any positive semidefinite matrix $G$ with ones on the diagonal,
\begin{equation}\label{equation_clip}
   \cR_\tau(G) \ge \min\{D(G),(\tau-1)^2\},\ 
   \text{ where }
D(G) = \tr(G-I)^2 = 2\sum_{i<j}|G_{ij}|^2.
\end{equation}
Indeed, if no eigenvalue exceeds $\tau$, the remainder equals $D(G)$. 
Otherwise, one eigenvalue contributes at least $\psi_\tau(\tau) = (\tau-1)^2$, 
and all other contributions are nonnegative. 
Also, for disjoint principal blocks $G_1,\ldots,G_l$ of $G$, we have
\begin{equation}\label{equation_cutR}
   \cR_\tau(G) \ge \sum_{j=1}^l\cR_\tau(G_j).
\end{equation}
To see this, we diagonalize each block in its coordinate subspace and 
complete it to an orthonormal basis. Its diagonal entries 
are averages of the eigenvalues of~$G$. 
Scalar Jensen's inequality for the convex function $\psi_\tau$ then
gives~\eqref{equation_cutR}; unused coordinates contribute nonnegative terms.

%%%%%%%%%%%%%%%%%%%%%%%%%%%%%%%%%
\begin{lemma}\label{Lemma_blocks}
If $114\epsilon\le(\tau-1)^2$, then
\begin{equation}\label{equation_blockremainder}
   \cR_\tau(H) \ge \frac{114\epsilon}{120}(n-119)-\frac{57}{2000}L.
\end{equation}
\end{lemma}
\begin{proof}
Inside a consecutive $120$-point block $J$, we sum the local inequality over 
its $114$ consecutive $7$-point span. By~\eqref{equation_wdiagonalsum}, 
any pair receives a total coefficient of at most two. Therefore
\begin{align*}
   D(H_J)+P_J \ge 114\epsilon,
 \ \text{ with }
 P_J = p\sum_{\substack{\text{$7$-point span }W\subset J\\
     W\ \text{consecutive}}}\operatorname{span}(W) \ge 0.
\end{align*}
The inequality~\eqref{equation_clip} and the assumed cap imply
$\cR_\tau(H_J)+P_J \ge 114\epsilon$.

For $n \ge 120$, we sum over all $n-119$ consecutive blocks $J$. 
Splitting their starting indices into the $120$ residue classes 
modulo $120$ gives disjoint blocks in each class. 
Hence~\eqref{equation_cutR} implies
\begin{align*}
   120\cR_\tau(H) \ge 114\epsilon(n-119)-\sum_JP_J.
\end{align*}
Each $7$-point span belongs to at most $114$ such blocks. 
Summing all $7$-point spans counts each elementary gap at most six times. 
Therefore, $\sum_JP_J\le114\cdot6pL$. 
Then, dividing by~$120$, we obtain~\eqref{equation_blockremainder}. 
For $n<120$, its right-side is nonpositive, so nonnegativity of the remainder suffices.
\end{proof}

For $\epsilon = 11/250$, the coefficient of $n-119$ is $\alpha = 209/5000$. 
With the selection specified in part~\\hypertarget{PartA}{}, 
the parameter $\tau$ from~\eqref{equation_tau} satisfies $\tau>13/4$, 
so that
\begin{equation}\label{equation_capcheck}
   114\epsilon = \frac{627}{125}<\frac{81}{16}<(\tau-1)^2.
\end{equation}
This proves~\eqref{equation_roadmapR}, our target in part~(\hyperlink{CaseB}{\texttt{B}}). 
Combining it with~\eqref{equation_roadmap} gives
\begin{equation}\label{equation_finitefinal}
   E_K(Z)+dL+119\alpha \ge (1+\alpha)N+\frac{\tau^2}{4}B.
\end{equation}

%%%%%%%%%%%%%%%%%%%%%%%%%%%%%%%%%%%%%%%%%%%%%%%%%%%%%%%%%%%%%%%
\section{The analytic input and completion of the proof}\label{Section_completion}
Next, we use the following precise fixed-test form of 
the unconditional Montgomery theorem~\cite[Lemma~5]{BGSTB2024}, 
which is also stated in~\cite[Lemma~3.1]{Lam2026}. If $q$ is real, even, integrable, 
supported in $[-1,1]$, and Lipschitz at zero, then
\begin{multline}\label{equation_PC}
   \sum_{0<\gamma,\gamma'\le T}
     \widehat q\left(\frac{i(\rho-\rho')\log T}{2\pi}\right)
               \frac4{4-(\rho-\rho')^2}\\
    = \frac{T\log T}{2\pi}
     \left(q(0)+\int_{-1}^1|v|q(v)\dd v
           + O_q((\log T)^{-1/2})\right).
\end{multline}
Both zero sums include multiplicity. The correction to the pointwise theorem 
in~\cite[Section~3]{BGSTB2026} explicitly leaves this integrated assertion unchanged.

For a fixed smooth even function $\eta\in C_c^\infty((-1/2,1/2))$ 
with $\int\eta^2 = 1$, put $f = \eta^2$, $K = \widehat f$, \mbox{and $Q = f*f$}. 
Then $\widehat Q = K^2$ and $\widehat{Q''}(z) = -4\pi^2z^2K(z)^2$. 
At $z = i(\rho-\rho')\log T/(2\pi)$,
\begin{equation*}%\label{equation_remove}
   K(z)^2 = \left(\widehat Q(z)-\frac{\widehat{Q''}(z)}{4(\log T)^2}\right)
                         \frac4{4-(\rho-\rho')^2}.
\end{equation*}
We apply formula~\eqref{equation_PC} separately to the two fixed smooth functions $Q$ and $Q''$. 
The second contributes $O_\eta(N(T)/(\log T)^2)$,
while the leading constant of the first is
\begin{align*}
   Q(0)+\int_{-1}^1|v|Q(v)\dd v = \cC(f).
\end{align*}
Evenness of $f$ changes $|u+v|$ to $|u-v|$ in the convolution integral. 
This proves~\eqref{equation_momentlimit} without 
any uniformity assumption for a changing test function.

Now we choose even smooth cutoffs $0\le\chi_\delta\le1$, 
supported inside $I$, and tending to one in its interior. 
We let
\begin{align*}
   \eta_\delta = \frac{\chi_\delta\sqrt{f_*}}{(\int\chi_\delta^2f_*)^{1/2}},
   \qquad f_\delta = \eta_\delta^2,\qquad K_\delta = \widehat f_\delta.
\end{align*}
The positivity of $f_*$ makes these admissible, and imply
\begin{align*}
   f_\delta\to f_*\ \text{in }L^1\cap L^2,\quad
   \cC(f_\delta)\to C,\ \text{ and }\
   \sup_{x\in\RR}|K_\delta(x)-K_*(x)|\le\|f_\delta-f_*\|_1\to0.
\end{align*}
Both kernels have an absolute value of at most one on $\RR$. 
Because the weights total twelve, the local inequality persists 
for $K_\delta$ with the smaller floor
\begin{align*}
   \epsilon_\delta = \frac{11}{250}-24\|f_\delta-f_*\|_1>0
\end{align*}
for small fixed $\delta$, and with the same gap penalty. 
As stated in part~(\hyperlink{CaseC}{\texttt{C}}), we take
\begin{align*}
   \alpha_\delta = \frac{114\epsilon_\delta}{120}\ \ \text{ and }\ \
   \tau_\delta = 2+\sqrt{2-2\alpha_\delta}.
\end{align*}
The cap condition remains valid, that is, 
its left-side decreases and its right-side increases relative 
to~\eqref{equation_capcheck}.

We apply the finite proof of~\eqref{equation_finitefinal}, 
with these parameters, to $Z_T$ from~\eqref{equation_ZT}. 
The functional equation sends $\rho$ to $1-\bar\rho$, 
which becomes conjugation under this normalization. 
Its simple real locations are in $[-L_T,0)$, where
\begin{align*}
   L_T = \frac{T\log T}{2\pi}\ \text{ and } L_T/N(T)\to 1
\end{align*}
by the Riemann--von Mangoldt formula~\cite[Section~2.2]{BGSTB2026}. 
The finite $B$ is exactly~\eqref{equation_badmass}. 
Dividing~\eqref{equation_finitefinal} by $N(T)$, 
first letting $T\to\infty$ for each fixed $\delta$, and then letting $\delta\to 0$, gives
\begin{equation}\label{equation_U}
   \liminf_{T\to\infty}\frac{\un(T)}{N(T)}
      \ge 1-\frac{4(C-1-\alpha+d)}{\tau^2}
    = 1-\frac{4(C-10133/10000)}{(2+\sqrt{4791}/50)^2}.
\end{equation}
Since we know from~\eqref{equation_C} that 
$C \approx 1.33133$, we find that the constant on the right-side of~\eqref{equation_U}
is approximately $0.88893 >8/9$,
which concludes the proof of Theorem~\ref{Theorem_main}.

\subsection*{Acknowledgements}
All calculations and figures were made in \texttt{SageMath}~\cite{SAGE}.

%%%%%%%%%%%%%%%%%%%%%%%%%%%%%%%%%%%%%%%%%%%%%%%%%%%%%%%%%%%%%%%%%%%
\bibliographystyle{plainurl}% shows urls
% \bibliography{Zeta89}% common bib file

% \end{document}

\vspace*{16mm}
\appendix
\section{\texorpdfstring{Algorithm that verifies the seven-point inequality on $\Omega'$}{Algorithm that verifies the seven-point inequality on Omega'}
}\label{Appendix}

%%%%%%%%%%%%%%%%%%%%%%%%%%%%%%%%%%%%
\vspace*{2mm}
\noindent
% \textbf{The algorithm.}
The goal is to prove $S(\bg)\ge S_0$ for $\bg\in\Omega'$, where $S_0=11/250=0.044$,
\begin{align*}
S(\mathbf g)=
\sum_{0\le i<j\le6}w_{ij}q(y_j-y_i)
+\frac1{200}\sum_{i=1}^6g_i,
\end{align*}
$q(x)=K(x)^2$, and $y_0=0$, $y_j=g_1+\cdots+ g_j$.

\vspace*{2mm}

\texttt{Step 1.}
Start with the closed box
$B_0=\left[743/1000, 5085/1000\right]^6$
and set the constraint to check 
$\sum_{i=1}^6g_i\le\frac{44}{5}$.

\vspace*{1mm}
\texttt{Step 2.}
For a box
$B=\prod_{i=1}^6[\ell_i,u_i]$,
discard it if 
$\sum_{i=1}^6\ell_i>\frac{44}{5}$.
Otherwise, tighten each upper endpoint by replacing it with
$u_k\leftarrow
\min\left(u_k,\frac{44}{5}
-\sum_{r\ne k}\ell_r\right)$.
This replacing is repeated in a loop until no further endpoint changes.

\vspace*{1mm}
% \noindent
\texttt{Step 3.}
For every $i<j$, calculate the interval 
$J_{ij} = \left[\sum_{k=i+1}^{j}\ell_k,\, \sum_{k=i+1}^{j}u_k \right]$ that
contains all values of $y_j-y_i$ on $B$.
Use interval arithmetic to bound $q$ on
$I_{ij}$, with the aim to estimate
\begin{align*}
L_1(B) :=
\sum_{i<j}w_{ij}\inf_{x\in J_{ij}}q(x)
+\frac1{200}\sum_{i=1}^6\ell_i.
\end{align*}
If $L_1(B)\ge S_0$, the goal inequality is proved on $B$.
(Actually, we can replace the infimum in $L_1(B)$ with some lower bound $m(I_{ij}) \le q(x)$ for $x\in I_{ij}$, which can be obtained using
monotonicity where it is established, and otherwise use  evaluations on iterated subdivisions of  $I_{ij}$. 
If the resulting lower bound is too weak to prove 
$L_1(B) \ge 11/250$, proceed to the next step with the Taylor lower-bound test or subdivide the six-dimensional box further.)

\vspace*{1mm}
\texttt{Step 4.}
If the first lower bound is insufficient, that is, $L_1(B)<S_0$, use a Taylor
expansion of $q$ at the center of $B$, with rigorous
bounds for $q'$ and $q''$. 
Calculate
\begin{align*}
    L_2(B) := S(\bc)-\sum_{k=1}^6|D_k|r_k
              -20\sum_{i<j}w_{ij}r_{ij}^2,
\end{align*}
where $\bc=(c_1,\dots,c_6)$, $\br=(r_1,\dots,r_6)$, and
$s_{ij}$,  $r_{ij}$, and $D_k$ 
are defined by~\eqref{equation_center-radii}.
Then, if $L_2(B)\ge S_0$, the inequality is proved on $B$.

If neither lower bound is $\ge S_0$, divide the longest
side of $B$ at its midpoint and add both resulting boxes
to the list of boxes to be examined.

\vspace*{1mm}
\texttt{Step 5.}
Repeat \texttt{Steps 2--5} until the list is empty or a
predefined resource limit is reached.

\vspace*{2mm}

There is no preset number of bisections of boxes. 
A finite successful run completes the proof of inequality~\eqref{equation_seven} when the list of boxes that remain to be checked becomes empty.
However, a lightweight implementation of the algorithm
in \texttt{SageMath} proved successful.
We remark that in order to ensure that the final result is correct, the variables
in question are evaluated, by case,  within safeguard lower or upper enclosures and appropriate rounding.
The running time of the program is under a minute, after which it finishes with zero unchecked boxes and the following output:

\vspace*{2mm}

\begin{verbatim}
    # N_processed	252,035
    # N_discarded	0
    # N_split	    126,017
    # N_direct in test L1(B)	44,195
    # N_Taylor in test L2(B)	81,823
\end{verbatim}
\smallskip
\begin{verbatim}
    The smallest lower bound for S in test L1(B):
        9907927838470733/225179981368524800
    The smallest lower bound for S in test L2(B):
        12384901590203241383/281474976710656000000
\end{verbatim}
\medskip

Note that The reported counters satisfy
$44195+81823 = 126018$ accepted boxes, and
twice the number of split boxes plus the initial box totals the number of processed boxes, that is,  the analyzis verified the nodes of a completed binary subdivision tree with no boxes discarded).
% $2\cdot 126017 + 1 = 126018.

The approximations of the least lower bounds obtained with the
acceptance tests  \texttt{L1(B)} and \texttt{L2(B)} were both above the target $0.044$, respectively:
\begin{align*}
    9907927838470733/225179981368524800&\approx
0.044000038450379066\\
12384901590203241383/281474976710656000000&\approx
0.044000009290113131.
\end{align*}

\vfill
\vspace*{6mm}

\end{document}